\documentclass[aps,12pt,nofootinbib]{revtex4}

\usepackage{graphicx,graphics}
\usepackage{enumerate}
\usepackage{subfigure}
\usepackage{amsmath}
\usepackage{amsfonts}
\usepackage{amssymb}
\usepackage{amsthm}
\usepackage{psfrag}

\catcode`\@=11
\def\numberbysection{\@addtoreset{equation}{section}
\def\theequation{\arabic{section}.\arabic{equation}}}
\numberbysection

\def\a{\alpha}
\def\d{{\rm d}}
\def\ci{{\rm i}}

\def\e{\varepsilon}
\def\Z{\mathbb{Z}}

\begin{document}

\title{Euler tours and unicycles in the rotor-router model}

\author{V.S. Poghosyan$^1$ and V.B. Priezzhev$^2$}
\affiliation{
$^1$Institute for Informatics and Automation Problems\\ NAS of Armenia, 0014 Yerevan, Armenia\\
$^2$Bogoliubov Laboratory of Theoretical Physics,\\ Joint Institute for Nuclear Research, 141980 Dubna, Russia
}

\begin{abstract}
A recurrent state of the rotor-routing process on a finite sink-free graph can be represented by a unicycle that is a
connected spanning subgraph containing a unique directed cycle. We distinguish between short cycles of length 2 called "dimers"
and longer ones called "contours". Then the rotor-router walk performing an Euler tour on the graph generates
a sequence of dimers and contours which exhibits both random and regular properties. Imposing initial conditions
randomly chosen from the uniform distribution we calculate expected numbers of dimers and contours
and correlation between them at two successive moments of time in the sequence. On the other hand, we prove that the excess
of the number of contours over dimers is an invariant depending on planarity of the subgraph but not on initial
conditions. In addition, we analyze the mean-square displacement of the rotor-router walker in the recurrent
state.
\end{abstract}

\maketitle


\section{Introduction}

The rotor-router walk is the latter and most frequently used name of the model introduced independently in different
areas during the last two decades. The previous names "self-directing walk"  \cite{P96} and "Eulerian walkers" \cite{PDDK} reflected its connection with the theory of self-organized criticality \cite{BTW} and the Abelian sandpile model \cite{Dhar}. Cooper and Spencer \cite{CS} called the model "P-machine" after Propp who proposed the rotor mechanism as the way to derandomize the internal diffusion-limited aggregation. Later on, several theorems in this direction have been proved in \cite{LP05,LP07,LP08}. Holroyd and Propp \cite{HP} proved a closeness of expected values of many quantities for simple random and rotor-router walks. Applications of the model to multiprocessor systems can be found in \cite{RSW}. Recent works on the rotor-router walk address the questions on recurrence \cite{AngelHol,HussSava}, escape rates \cite{FGLP} and transitivity of the rotor-routing action \cite{CCG}.

The connection between the Abelian sandpiles, Euler circuits and the rotor-router model observed in the original paper \cite{PDDK} was the subject of the rigorous mathematical survey \cite{HLMPPW}. An essential idea highlighted in the
survey is the consideration of the rotor-routing action of the sandpile group on spanning trees in parallel with
rotor-routing on unicycles. The rotor-router walk started from an arbitrary rotor configuration on a finite sink-free
directed graph $G$ enters after a finite number of steps into an Euler circuit ( Euler tour) and remains there forever.
The length of the circuit is the number of edges of the digraph.
Each recurrent rotor state can be represented by a connected
spanning subgraph $\rho\subset G$ which contains as many edges
as vertices and contains a unique directed cycle \cite{CRST,RetProb,LoopingConst}.
The dynamics of the rotor-router walk requires the location of the walker at a vertex $v\in \rho$ belonging to the cycle.
The pair $(\rho,v)$ is called unicycle (see Section II for precise
definition). Thus, the walk passes the periodic sequence of unicycles.

\begin{figure}[h!t]
  \centering
  \includegraphics[width=120mm]{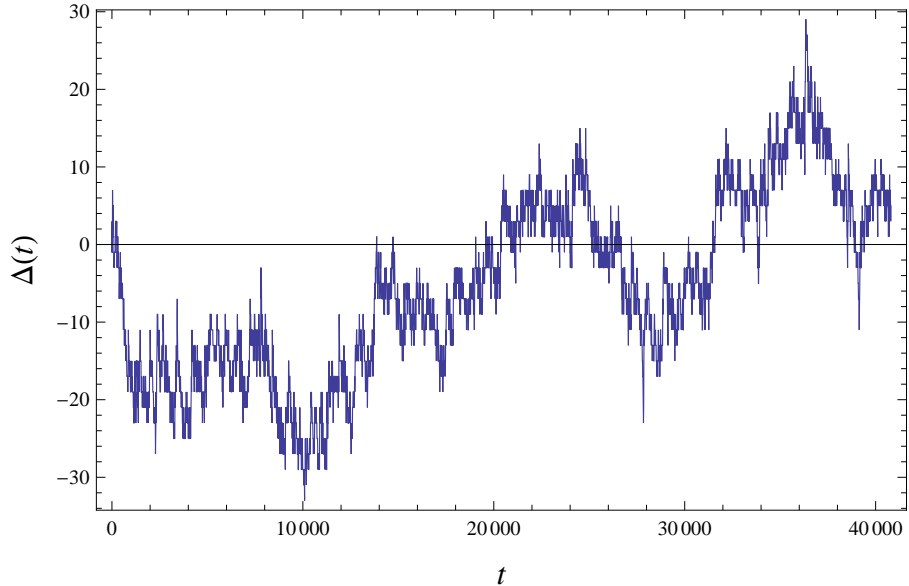}
  \caption{"Displacement" $\Delta(t)$ for the Euler tour of 40000 steps.}
  \label{Delta1}
\end{figure}

A shortest cycle in the unicycle is the two-step path from a given vertex to one of nearest neighbors and back.
We call the cycles of length 2 "dimers" by analogy with lattice dimers covering two neighboring vertices. Longer
cycles involve more than two vertices and form directed contours. The Euler tour
passes sequentially unicycles containing cycles of different length. The order in which dimers and contours alternate
depends on the structure of the initial unicycle. Ascribing $+1$ to each step producing a contour and $-1$ to a dimer,
we obtain for a "displacement" $\Delta(t)$ after $t$ time-steps the picture (Fig.\ref{Delta1}) resembling the symmetric random walk.
Nevertheless, the process actually is neither completely symmetric nor completely random.
It is the aim of the present paper to investigate statistical properties
of unicycles as they appear in course of the Euler tour.

We will see that the events "dimer" and "contour"  correlate along the sequence and an excess of the number of dimers
over contours is an invariant characterizing topology of the surface where the rotor-router walk occurs. Specifically, in the limit of large square lattice with periodic boundary conditions we find the expected number of dimers in the Euler tour and an analytical expression for the correlations dimer-dimer and contour-contour at two successive moments of time in the circuit. We consider a closed loop encircling a plane domain and prove that the rotor-router walk
passed each directed edge of the domain contains the number of dimers exceeding that of contours exactly by $1$.
This property does not hold for surfaces of the non-zero genus.

In addition to statistics of unicycles, we consider the mean-square displacement of the rotor-router walker in the recurrent state and argue that it yields to the diffusion law with the diffusion coefficient depending on dynamic
rules and boundary conditions.

\section{The model}

Consider a directed graph (digraph) $G=(V,E)$ with the vertex set $V=V(G)$ and the set of directed edges $E=E(G)$ without self-loops and multiple edges.
If for each edge directed from $v$ to $w$, there exists an edge directed from $w$ to $v$, graph $G$ is bidirected.
The bidirected graph can be obtained by replacing each edge of an undirected graph with a pair of directed edges, one in each direction.

A subgraph $G\,'$ of a digraph $G$ is a digraph with vertex set $V(G\,') = V(G)$ and edge set $E(G\,')$ being a subset of $E(G)$, i.e. $E(G\,') \subseteq E(G)$. In this case we write $G\,' \subseteq G$. If $E(G\,')$ contains no outgoing
edges from a fixed vertex, that vertex is a sink. The oriented tree with sink $v$ is a digraph, which is acyclic and whose every non-sink vertex $w \neq v$ has only one outgoing edge.
If the subgraph of $G$ is a tree with sink $v$ then it is called a spanning tree of $G$ with root $v$.
A connected subgraph of an oriented graph $G$, in which every vertex has one outgoing edge, is called unicycle.
The unicycle contains exactly one directed cycle.

An  Euler circuit (or Euler tour) in a directed graph is a path that visits each directed edge exactly once.
If such a path exists, the graph is called Eulerian digraph.
A digraph is strongly connected if for any two distinct vertices $v$, $w$ there are directed paths from $v$ to $w$ and from $w$ to $v$.
A strongly connected digraph $G=(V,E)$ is Eulerian if and only if for each vertex $v \in V$ the in-degree and out-degree of $v$ are equal.
In particular, the one-component bidirected graph is Eulerian.
We call $G$ an Eulerian digraph with sink if it is obtained from an Eulerian digraph by deleting all the outgoing edges from one vertex.
The subset of sites of $G$ connected with the sink forms an open boundary.

The rotor-router model is defined as follows. Consider an arbitrary digraph $G = (V,E)$.
Denote the number of outgoing edges from the vertex $v \in V$ by $d_v$. The total number of edges of $G$ is $|E| = \sum_{v\in V} d_v$. Each vertex $v$ of $G$ is associated with an arrow, which is directed along one of the outgoing edges from $v$. The arrow directions at the vertex $v$ are specified by an integer variable $\a_v$, which takes the values $0 \leq \a_v \leq d_v - 1$. The set $\{ \a_v|\; v\in V,\; 0 \leq \a_v \leq d_v - 1\}$ defines the rotor configuration (the medium). Starting with an arbitrary rotor configuration one drops a chip to a vertex of $G$ chosen at random. At each time step the chip arriving at a vertex $v$, first changes the arrow direction from $\a_v$ to $ (\a_v + 1) \bmod d_v$\,, and then moves one step along the new arrow direction from $v$ to the corresponding neighboring vertex. The chip reaching the sink leaves the system. Then, the new chip is dropped to a site of $G$ chosen at random.

In the absence of sinks the motion of the walker does not stop. The rotor configuration $\rho$ can be considered as a subgraph of $G$ $(\rho \subset G)$ with the set of vertices $V(\rho) = V$ and the set of edges $E(\rho) \subset E$ obtained from the arrows. The state of the system (single walker + medium) at any moment of time is given by the pair $(\rho, v)$ of the rotor configuration $\rho$ and the position of the chip $v \in V$. According to arguments
in \cite{PDDK}, the rotor-router walk started from an arbitrary initial state $(\rho, v)$ passes transient states and
enters into a recurrent state, continuing the motion in the limiting cycle which is the Eulerian circuit of the graph.
The basic results about the rotor-router model on the Eulerian graphs can be summarized as two propositions.

{\it Proposition 1} [\cite{HLMPPW}, Theorem 3.8] Let $G$ be a strongly connected digraph. Then $(\rho, v)$ is a
recurrent single-chip-rotor state on $G$ if and only if it is a unicycle.

The rotor states that are not unicycles are transient.
In contrast to recurrent states, they appear at the initial stage of evolution
up to the moment when the system enters into the Eulerian tour.

{\it Proposition 2} [\cite{HLMPPW}, Lemma 4.9] Let $G$ be an Eulerian digraph with $m$ edges. Let $(\rho, v)$ be
a unicycle in $G$. If one iterates the rotor-router operation $m$ times starting from $(\rho, v)$, the chip traverses
an Euler tour of $G$, each rotor makes one full turn, and the state of the system returns to $(\rho, v)$.

\section{The unicycles on torus}

Below, we specify the structure of graph $G$ as the square $N\times M$ lattice with periodic boundary
conditions (torus). Then the number of outgoing edges is $4$ for all vertices $v \in G$. We consider two ways
of labeling of four directions of the rotor $\alpha_v=0,1,2,3$ corresponding to the clockwise and cross order
of routing. For the clockwise routing, we put $\{ 0\equiv North,1\equiv East,2\equiv South, 3\equiv West\}$,
for the cross one $\{ 0\equiv North,1\equiv South,2\equiv East, 3\equiv West\}$. Sometimes, we will use the
notation $\alpha_v$ for the edge outgoing from $v$ in direction $\alpha_v$.

By {\it Proposition 1}, each recurrent state $(\rho, v)$ of the rotor-routing process is a unicycle on $G$.
Replacing all arrows of the rotor configuration $\rho$ by directed bonds, we obtain a cycle-rooted spanning tree
of $G$. It consists of a single cycle of length $s$ (i.e. $s$ directed bonds connecting $s$ vertices) and the
spanning cycle-free subgraph whose edges are directed towards the cycle. If $s=2$, the cycle is "dimer", if
$s > 2$, the cycle is "contour" oriented clockwise or anticlockwise.

By {\it Proposition 2}, the walker started from unicycle $(\rho, v)$ traverses an Euler tour which has the length
$m=4MN$ in our case of $M\times N$ torus.
The questions arise how many of $m$ unicycles passed during the Euler tour contain a dimer (contour)?
What is the probability that the unicycle obtained on $t$-th step of the Euler tour contains a dimer (contour)?

Consider the recurrent state $(\rho, v)$ and define a random variable $X_{\alpha_v}$ as
\begin{equation}
X_{\alpha_v} = \left\{ \begin{array}{cc}
                1 & \text{if the cycle containing edge $\alpha_v$ is dimer } \\
                0 & \text{if the cycle containing edge $\alpha_v$ is contour}
              \end{array}
\right. .
\end{equation}
We are interested in the average $<X_{\alpha_v}>$ over all possible uniformly distributed recurrent states and
over all directions $\alpha_v=0,1,2,3$.

We write the unicycle $(\rho, v)$ as $(T_v,\alpha_v,v)$ separating off the edge $\alpha_v$ and the spanning tree
$T_v$ obtained from edges outgoing from vertices of the set $V\smallsetminus v$.
By the definition of Euler tour, all $m$ unicycles following the initial unicycle $(T_v,\alpha_v,v)$ are different.
If two Euler tours have a common element, they coincide.

Now, we fix a vertex $w\in V$ and its outgoing edge $\alpha_w$. If one scans
over all possible initial spanning trees $T_v$, then the trees $T_w$ also scan over all possible configurations.
So, the uniform distribution of  $T_v$ induces the uniform distribution of  $T_w$. Therefore, $<X_{\alpha_w}>$ is
the probability that the edge $\alpha_w$ taken uniformly with $\alpha_w=0,1,2,3$ and added to the uniformly
distributed spanning trees creates a dimer. Due to
the translation invariance, this average does not depend on the position of the initial vertex $v$.

To make these arguments more explicit, consider all possible unicycles
\begin{equation}
(\rho,w) \equiv (T_w,\alpha_w,w)
\end{equation}
for fixed vertex $w$ and arrow $\alpha_w$.
First, we take $\alpha_w=0$, choosing the arrow at $w$ directed {\it North}.
The set of unicycles $(T_w,0,w)$ can be divided into two subsets $(T_w,0,w)_d$ and $(T_w,0,w)_c$ where the first subset corresponds to
spanning trees $T_w$ containing a selected bond incident to $w$ from above.
The tree $T_w$ has the root in $w$, so this bond is directed down.
The selected bond and arrow $\alpha_w=0$ form together a vertical dimer with the lower end in $w$. In the subset $(T_w,0,w)_c$,
the place of the selected bond in each tree $T_w$ is empty, so the arrow $\alpha_w=0$ belongs to a contour.
Considering similar subsets for other directions $\alpha_w=1,2,3$ with selected bonds of the trees $T_w$ incident to $w$
from right, down and left, we can write the average probability to find a dimer incident to $w$ as
\begin{equation}
P(d)= \frac{1}{4|T|}\sum_{\alpha_w}(T_w,\alpha_w,w)_d,
\end{equation}
where summation is over all $\alpha_w=0,1,2,3$ and
$|T|$ is the total number of non-rooted spanning trees. Now, let us take the sum over all $w$ in the numerator and denominator
using the uniformity of  vertices of the torus.
Then, the numerator will be the doubled number of edges of the spanning tree $|E_T|$ multiplied by $|T|$ because each edge is taken in two directions.
The denominator will be $4MN|T|$.
The number of edges of the torus is $|E| = 2 MN$, the number of edges of the spanning tree $|E_T| = MN - 1$.
Therefore, the probability of a dimer $P(d)$ is
\begin{equation}
P(d) \equiv <X_{\alpha_w}> = \frac{|E_T|}{|E|} =
 \frac{1}{2} - \frac{1}{2 MN}\,.
\label{Probdim}
\end{equation}
The probability of a contour is $P(c)=1-P(d)$.

In the limit  $M\rightarrow \infty$, $N\rightarrow \infty$, we obtain $P(c)=P(d)=1/2$.
In spite of this simple symmetric result, the distribution of the random value $X_{\alpha_v}$ is not trivial.
We will return to this question in the next section.
Now consider the correlations dimer-dimer and dimer-contour at two successive moments of time in the Euler tour.

In the case of clockwise routing where the directions of each arrow alternate $North$--$East$--$South$--$West$,
two successive directions of an arrow at a fixed vertex always form the angle $90^{\circ}$.
For the cross routing ($North$--$South$--$East$--$West$), the rotations $North$--$South$ and $East$--$West$ form the angle $180^{\circ}$,
whereas the rotations $South$--$East$ and $West$--$North$ form the angle $90^{\circ}$.

The correlations we are going to determine have the following origin. Consider for example a particle arriving to vertex $v$
from above at the time step $t$. If the arrow at $v$ is directed {\it North} in the preceding moment of time, a vertical dimer
is created with the lower vertex in $v$. If one uses the clockwise dynamics, the next step is the rotation of the arrow at $v$
to {\it East}. Assume that there is an arrow at time $t$ directed to $v$ from right to left. Then, the horizontal dimer
is created at the time step $t+1$ with the left vertex in $v$. The probability to get two dimers at moments $t$ and $t+1$
is the correlation $P(d,d)$.

The arguments used for the derivation $P(d)$ and $P(c)$ show that the correlations $P(d,d)$, $P(d,c)$, $P(c,d)$ and $P(c,c)$ at two successive time-steps can be related with the probability to find two adjacent edges of the square lattice occupied (or not occupied) by bonds of the spanning tree $T$.

Specifically in the considered example, we must enumerate unicycles $(T_v,0,v)$ with the vertical dimer having the lower end in $v$. The spanning tree $T_v$ of the unicycle $(T_v,0,v)$ has the root $v$ and two fixed bonds, $b_1$ directed to $v$ from above
and $b_2$ directed to $v$ from the right. The presence of the root in $v$ implies that all bonds of the tree $T_v$ are globally
oriented towards the vertex $v$.

The enumeration of spanning trees $T_v$ obeying the above conditions can be performed in three steps. First, we consider non-oriented spanning tree $T$ with the selected non-oriented bonds $\bar{b_1}$ and $\bar{b_2}$ on the places of $b_1$ and $b_2$.
Second, we put the root in $v$, giving the necessary orientation to bonds $b_1$ and $b_2$ and supplying  other bonds with the global orientation towards $v$. Third, we use the Kirhhoff theorem according to which the number of spanning trees does not depend on the location of the root. This allows us to shift the root to infinity and restore the translation invariance in the limit
$M\rightarrow \infty$ and $N\rightarrow \infty$.

The alternative way of calculations would consist in fixing the oriented bonds $b_1$ and $b_2$ and the location of the root in $v$. However, the Kirhhoff theorem does not allow the translation of the root in the presence of oriented bonds. Then, the lack
of translation invariance makes all calculations much more difficult.

We fix a vertex $i_0 \in V$ and its two neighbors on the square lattice $i_1$ and $i_2$. Then
$e_1 = \{i_0,i_1\}$ and $e_2 = \{i_0,i_2\}$ are adjacent edges.

Define the probabilities
\begin{eqnarray}
P(++) &=& \mathrm{Prob}(e_1 \in T, e_2 \in T),\\
P(-+) &=& \mathrm{Prob}(e_1 \notin T, e_2 \in T),\\
P(+-) &=& \mathrm{Prob}(e_1 \in T, e_2 \notin T),\\
P(--) &=& \mathrm{Prob}(e_1 \notin T, e_2 \notin T).
\end{eqnarray}
Obviously, $P(++) + P(-+) + P(+-) + P(--) = 1$ and $P(-+) = P(+-)$ due to symmetry.

The calculation of probabilities of fixed spanning tree configurations is a standard procedure, which uses
the Green functions and so called defect matrices (see e.g. \cite{MajDhar,Pri94,PirouxRu,P12,PogPri}). In our case, it gives
\begin{eqnarray}
P(++) &=& \lim_{\e \rightarrow \infty} \frac{\det(I + B_1 G)}{\e^2},\\
P(+-) &=& P(-+) = \lim_{\e \rightarrow \infty} \frac{\det(I + B_2 G)}{\e},\\
P(--) &=& \det(I + B_3 G),
\end{eqnarray}
where the matrices $I$, $G$ are
\begin{equation}
I =
\left(
\begin{array}{ccc}
 1 & 0 & 0 \\
 0 & 1 & 0 \\
 0 & 0 & 1
\end{array}
\right)
,\quad
G =
\left(
\begin{array}{ccc}
 G_{i_0,i_0} & G_{i_0,i_1} & G_{i_0,i_2} \\
 G_{i_1,i_0} & G_{i_1,i_1} & G_{i_1,i_2} \\
 G_{i_2,i_0} & G_{i_2,i_1} & G_{i_2,i_2} \\
\end{array}
\right),
\end{equation}
and the defect matrices  $B_1$, $B_2$ and $B_3$ are
\begin{equation}
B_1 = \left(
\begin{array}{ccc}
 2\,\e & -\e & -\e \\
 -\e & \e & 0 \\
 -\e & 0 & \e \\
\end{array}
\right)
,\quad
B_2 = \left(
\begin{array}{ccc}
 \e-1 & -\e &  1 \\
  -\e &  \e &  0 \\
    1 &   0 & -1 \\
\end{array}
\right)
,\quad
B_3 = \left(
\begin{array}{ccc}
 -2 &  1 &  1 \\
  1 & -1 &  0 \\
  1 &  0 & -1 \\
\end{array}
\right).
\end{equation}

Defect matrices define the locations of bonds $e_1$ and $e_2$ which form angles  $90^{\circ}$ or  $180^{\circ}$.
In the first case we add index $a$ to the notations of probabilities, and index $b$ for the second case.
Using the explicit values for the Green functions given in Appendix, we obtain
in the limit $M\rightarrow \infty$ and $N\rightarrow \infty$
\begin{eqnarray}
P_a(++) &=& P_a(--) = \frac{1}{\pi} - \frac{1}{\pi^2},\\
P_a(-+) &=& P_a(+-) = \frac{1}{2} - \frac{1}{\pi} + \frac{1}{\pi^2}
\end{eqnarray}
for the case (a), and
\begin{eqnarray}
P_b(++) &=& P_b(--) = \frac{2}{\pi} - \frac{4}{\pi^2},\\
P_b(-+) &=& P_b(+-) = \frac{1}{2} - \frac{2}{\pi} + \frac{4}{\pi^2}
\end{eqnarray}
for the case (b).

Then, for the correlations dimer-dimer and dimer-contour at two successive moments of time in the Euler tour we have
\begin{eqnarray}
P(c,c)=P(d,d) = P_a(++)
&=& \frac{1}{\pi} - \frac{1}{\pi^2},\\
P(c,d)=P(d,c) = P_a(-+)
&=& \frac{1}{2} - \frac{1}{\pi} + \frac{1}{\pi^2}
\end{eqnarray}
in the case of the clockwise routing, and
\begin{eqnarray}
P(c,c)=P(d,d) = \frac{P_a(++) + P_b(++)}{2}
&=& \frac{3}{2\pi} - \frac{5}{2\pi^2},\nonumber\\\\
P(c,d)=P(d,c) = \frac{P_a(-+) + P_b(-+)}{2}
&=& \frac{1}{2} - \frac{3}{2\pi} + \frac{5}{2\pi^2},\nonumber\\
\end{eqnarray}
in the case of cross routing.

\section{The balance between dimers and contours}

Consider a part of the Euler tour $E(\rho_1,v_1|\rho_2,v_2)$ as a sequence of unicycles with the first element
$(\rho_1,v_1)$ and the last element $(\rho_2,v_2)$. The whole Euler tour in these notations is $E(\rho,v|\rho,v)$ and
the last unicycle $(\rho,v)$ is not included into the sequence.
We define a random value $\Delta(\rho_1,v_1|\rho_2,v_2)$ as
\begin{equation}
\Delta(\rho_1,v_1|\rho_2,v_2) = \mathrm{\#\,contours - \#\,dimers\; in\;} E(\rho_1,v_1|\rho_2,v_2)\,.
\end{equation}

If the cycle $C$ in unicycle $(\rho,v)$ is a contour oriented clockwise, we denote by $(\bar{\rho},v)$ the unicycle which differs from $(\rho,v)$ only by the counter-clockwise orientation of the contour. The following proposition for
the clockwise routing has been announced in \cite{PPS} and formulated in \cite{HLMPPW} as {\it Corollary }:

Let $G$ be a bidirected planar graph and let $(\rho,v)$ be a unicycle with the cycle $C$ oriented clockwise.
After the rotor-router operation is iterated some number of times, each rotor internal to $C$ has performed a full rotation, each rotor external to $C$ has not moved, and each rotor on $C$ has performed a partial rotation so
that the cycle is counter-clockwise  $\bar{C}$.

Below, we prove that $\Delta(\rho,v|\bar{\rho},v)=-1$ for any $\rho$ and $v \in V$ if the subgraph surrounded by $C$
is planar and the walker moves according to the clockwise routing.
It is important to note, that the clockwise routing is crucial both for the {\it Corollary} and the identity
$\Delta(\rho, v|\bar{\rho}, v) = -1$.

Using the method of induction, we start with the case of minimal $C$ when the contour is an elementary square $C_1$
of area 1. In this case, the walker makes four steps: starts with the clockwise contour $C_1$, produces sequentially
three dimers and ends by counter-clockwise contour $\bar{C_1}$. Denoting this sequence as $c,d,d,d,c$ we see that
$\Delta(\rho,v|\bar{\rho},v)=-1$.

Consider now unicycle $(\rho_s, v)$ containing a clockwise contour $C_s$ of area $s > 1$. We assume that
$\Delta(\rho_{s'},v|\bar{\rho_{s'}},v)=-1$ for all $1 < s' \leq s $ and prove  that $\Delta(\rho_{s+1},v|\bar{\rho}_{s+1},v)=-1$. Due to the symmetry of the square lattice, we can fix without loss
of generality the edge $\alpha_v$ with $\alpha_v=North$ at the left side of the contour $C_s$. Let $(x,y)$ be
coordinates of the vertex $v$. We consider four stages of the transformation of unicycles from $(\rho_{s+1},v)$ to
$(\bar{\rho}_{s+1},v)$.

{\it Stage 1.} The chip moves from $(x,y)$ to $(x+1,y)$ along the edge $\alpha_{(x,y)}=East$. If the obtained cycle is
dimer, Stage 1 is completed. Otherwise, the cycle is a contour $C_{s_1}$ of area $s_1 \leq s$. By {\it Corollary}, the unicycle with contour $C_{s_1}$ transforms after some number of steps into the unicycle with $\bar{C}_{s_1}$ and
by the assumption, $\Delta(\rho_{s_1},(x,y)|\bar{\rho}_{s_1},(x,y))=-1$. The contour $\bar{C}_{s_1}$ contains the
edge $\alpha_{(x+1,y)}=West$.

{\it Stage 2.} The chip moves from $(x+1,y)$ to $(x+1,y+1)$ along the edge $\alpha_{(x+1,y)}=North$. If the obtained cycle is dimer, Stage 2 is completed. Otherwise, the cycle is a contour $C_{s_2}$ of area $s_2 \leq s$.
By {\it Corollary}, the unicycle with contour $C_{s_2}$ transforms into one with $\bar{C}_{s_2}$ and
by the assumption, $\Delta(\rho_{s_2},(x+1,y)|\bar{\rho}_{s_2},(x+1,y))=-1$. The contour $\bar{C}_{s_2}$ contains the
edge $\alpha_{(x+1,y+1)}=South$.

{\it Stage 3.} The chip moves from $(x+1,y+1)$ to $(x,y+1)$ along the edge $\alpha_{(x+1,y+1)}=West$. If the obtained cycle is dimer, Stage 3 is completed. Otherwise, the cycle is a contour $C_{s_3}$ of area $s_3 \leq s$.
Again, the unicycle with contour $C_{s_3}$ transforms into one with $\bar{C}_{s_3}$ and
$\Delta(\rho_{s_3},(x+1,y+1)|\bar{\rho}_{s_3},(x+1,y+1))=-1$. The contour $\bar{C}_{s_3}$ contains the
edge $\alpha_{(x,y+1)}=East$.

{\it Stage 4.} The chip moves from $(x,y+1)$ to $(x,y)$ along the edge $\alpha_{(x,y+1)}=South$ and produce
the original unicycle but with the opposite orientation of the contour.

The description of evolution of unicycles from $(\rho,v)$ to $(\bar{\rho},v)$ shows that the number of dimers exceed that of contours by 1 during every stages 1,2,3 independently of whether the first or the second scenario of evolution is realized at each stage. Taking into account that the cycles of the first unicycle $(\rho,v)$ and the last one
$(\bar{\rho},v)$ are contours, we obtain
\begin{equation}
 \Delta(\rho,v|\bar{\rho},v)=-1\,.
 \label{theorem}
\end{equation}

{\it Remark}. \,\, According to {\it Corollary} the sum $s_1+s_2+s_3=s$\,.

\begin{figure}[h!t]
  \centering
  \includegraphics[width=120mm]{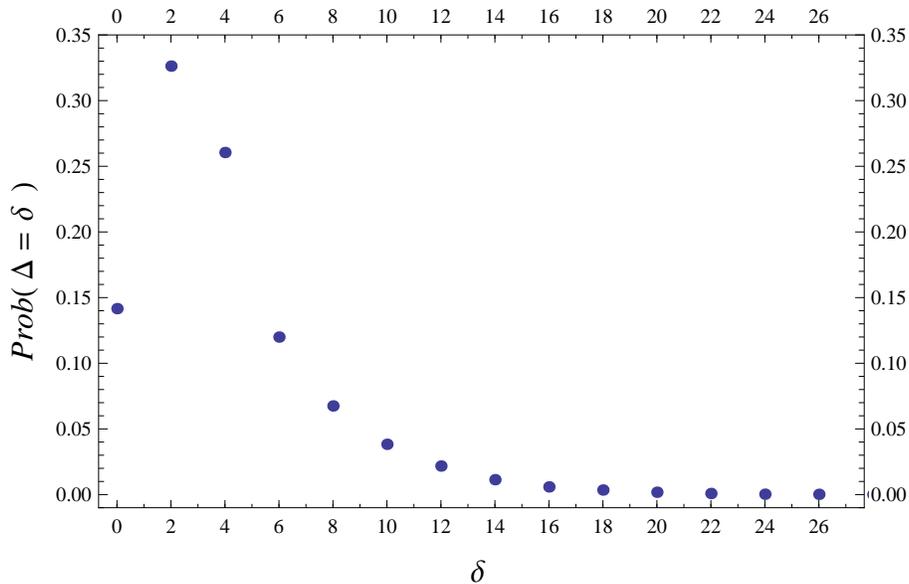}
  \caption{Probability that the number of contours exceeds the number of dimers by $\delta$ during
  the whole Euler tour on the torus $N=M=100$ for clockwise routing. The number of samples in simulation $10^6$.
  No negative $\delta$'s are detected.
   }
  \label{Fig-contour-distrib-clock}
\end{figure}
The result (\ref{theorem}) proven for the plane domain is in a drastic contrast with the rotor-router walk on surfaces
of the non-zero genus. In Fig.\ref{Fig-contour-distrib-clock} the function $\Delta(\rho,v|\rho,v)$ is shown for the
whole Euler tour on the torus for clockwise routing. We see that $\Delta(\rho,v|\rho,v)$ takes different values
depending on $(\rho,v)$ and all these values are non-negative. The average $\langle \Delta \rangle$ is known
from (\ref{Probdim}). Indeed, for the Euler tour of length $m=4MN$, we have the average number of dimers $2MN - 2$
and the average number of contours $2MN+2$. Therefore $\langle \Delta \rangle$=4 for any $M>1$ and $N>1$.

To consider a more general situation, we fix a unicycle $(\rho,v)$ with a clockwise contour $C$ which cuts out
a plane domain $A$ from the torus. According to {\it Corollary}, the contour $C$ in  $(\rho,v)$ will be converted
into the counter-clockwise contour $\bar{C}$ in $(\bar{\rho},v)$ after some number of steps of the Euler tour
started from $(\rho,v)$. The contour $\bar{C}$ is counter-clockwise with respect to $A$ and clockwise with respect to
the complement domain $A^c$ of genus 1. Now, we separate the whole Euler tour into two parts: from $(\rho,v)$ to
$(\bar{\rho},v)$ and from $(\bar{\rho},v)$ to $(\rho,v)$. From (\ref{theorem}) we have $\Delta(\rho,v|\bar{\rho},v)=-1$.
Therefore, to provide non-negativity of $\Delta(\rho,v|\rho,v)$, we should admit $\Delta(\bar{\rho},v|\rho,v)\geq 2$.
A reason for the excess of contours over dimers is the existence of many additional loops on the surface of
non-zero genus, in particular non-contractible loops on the torus. We are not able to prove an exact inequality
for $\Delta(\bar{\rho},v|\rho,v)$, so we formulate it as a conjecture:

{\it Conjecture} Let $\bar{C}$ be a contour clockwise with respect to the surface of genus 1 and $(\bar{\rho},v)$ is unicycle containing $\bar{C}$. Then, for the sequence of unicycles $(\bar{\rho},v),\dots,(\rho,v)$ in the Euler tour, the difference $\Delta(\bar{\rho},v|\rho,v)\geq 2$.

\begin{figure}[h!t]
  \centering
  \includegraphics[width=120mm]{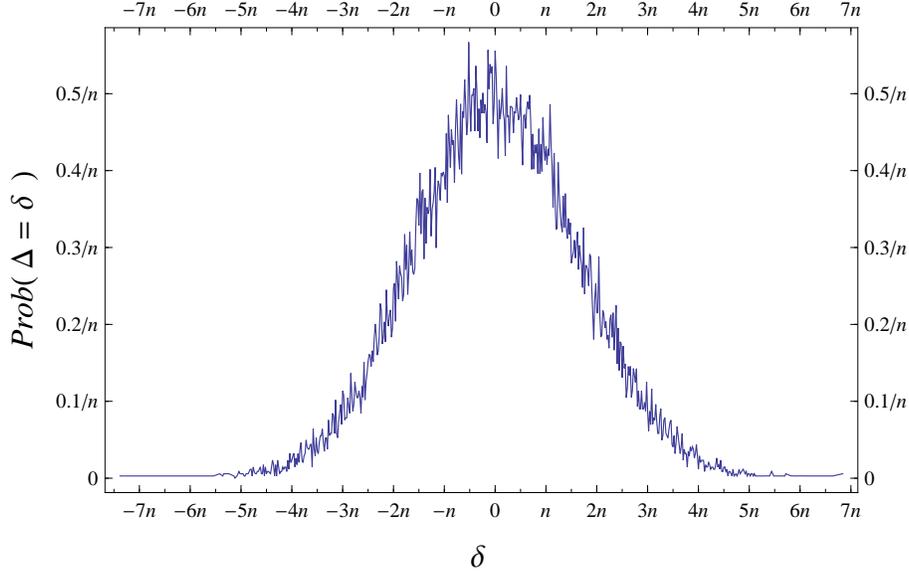}
  \caption{Probability that the number of contours exceeds the number of dimers by $\delta$ during
  the whole Euler tour on the torus $n \equiv N = M = 100$ for cross routing. The number of samples in simulation $10^6$.
  }
  \label{Fig-contour-distrib-cross}
\end{figure}

The conjectured inequality as well as (\ref{theorem}) get broken in the case of cross routing. Fig.\ref{Fig-contour-distrib-cross} shows $\Delta(\rho,v|\rho,v)$ for the Euler tour with the cross routing rules.
Instead of the strictly asymmetric distribution in Fig.\ref{Fig-contour-distrib-clock}, we have a Gaussian-like
distribution with the width corresponding to the diffusion law. To check the Gaussian nature of the distribution,
we calculated moments $m_2,m_3,m_4$, and estimated skewness and excess kurtosis. For the lattice size $n\equiv M = N = 100$ with
statistics of $10^6$ samples, we obtained $m_3/m_2^{3/2}=0.027594$ and $m_4/m_2^2-3=0.00594$. Nevertheless, the exact normality of the distribution in the limit $M\rightarrow \infty$ and $N\rightarrow \infty$ remains an unproved conjecture.

The average $\langle \Delta \rangle=4$ coincides with that for the clockwise routing because the probability
(\ref{Probdim}) does not depend on the order of routing.

\section{The diffusion of the walker}

Given the Euler tour of length $4 n^2$ on the torus $n \times n$, we can find the mean-square displacement
$<r(t)^2>$ after $t$ steps, where  $\vec{r}(t) = (x(t), y(t))$ and $ x(t), y(t)$ are coordinates of the walker
at time  $0 \leq t \leq 4n^2$. Fig.\ref{Fig-clock-torus} shows $<r(t)^2>$ for two periods of the Euler tour with
the clockwise routing. The interpolation of the function $<r(t)^2>$ in the interval $1 \ll t \ll n^2$ gives the linear dependence $<r(t)^2> \sim t $.

The obtained linear law is not surprising. The time dependence of mean square displacement cannot be slower than
$ct$, where $c$ is a constant. Indeed, by the definition of Euler tour, each vertex of the torus cannot be visited
more than 4 times. Therefore, the walker cannot stay in an area of radius $r$ longer than $4r^2$ time steps.

On the other hand, $<r(t)^2>$ cannot be faster than $kt$ where $k$ is another constant. It follows from the
{\it Corollary} that the walk is "loop-fiiling", i.e. the interior of a loop of radius $r$ is visited densely, so
that each rotor inside the loop makes a full rotation before the walker leaves the loop. Therefore, an advance of the walker at the distance of order $r$ takes $\sim r^2$ steps.

The exact value of the diffusion constant is unknown. The computer simulations show that it depends on the order of routing and we can estimate it as:
\begin{figure}[h!t]
  \centering
  \includegraphics[width=120mm]{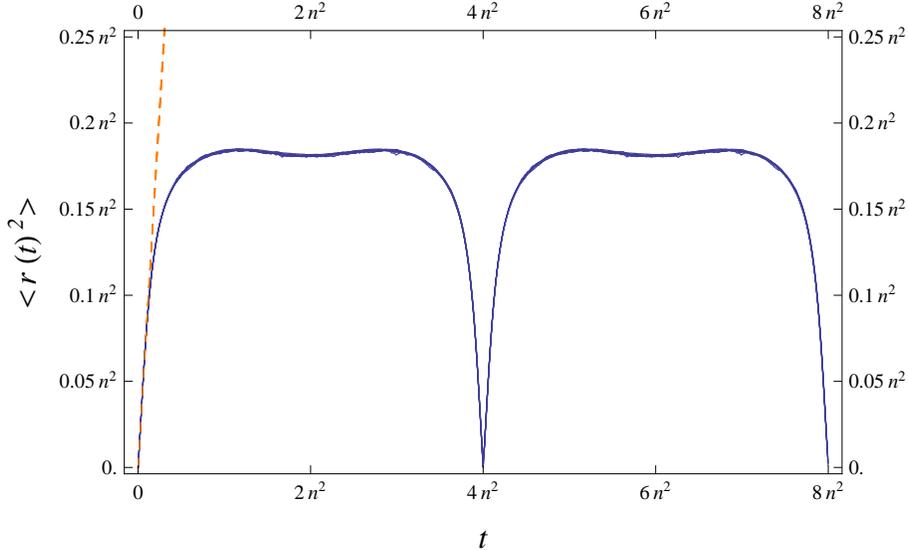}
  \caption{ The mean square displacement $<r(t)^2>$ for clockwise routing on the torus of size $n\times n$, $n=100$. The linear interpolation in the interval $0 \leq t \leq n$ gives $<r(t)^2> \simeq 0.83\,t $. The initial conditions are uniform spanning unicycles. The number of samples in simulation $10^6$. }
  \label{Fig-clock-torus}
\end{figure}
\begin{eqnarray}
<r(t)^2>&\simeq& 0.83\,t, \text{\,\,for clockwise routing\,,}\\
<r(t)^2>&\simeq& 1.32\,t, \text{\,\,for cross routing\,.}
\label{diffusion}
\end{eqnarray}

It is important to note, that the diffusion law (\ref{diffusion}) for the linear part of the Euler tour differs
from the subdiffusion law $<r(t)^2> \sim t^{2/3}$ obtained in \cite{PDDK} for the rotor walk in the infinite random media. The rigorous proof of the exponent $2/3$ is a challenging problem of the theory.

\section*{Acknowledgments}
This work was supported by the RFBR grants No. 12-01-00242a, 12-02-91333a,
the Heisenberg-Landau program, the DFG grant RI 317/16-1 and
State Committee of Science MES RA, in frame of the research project No. SCS 13-1B170.

\section{Appendix}

The translation invariant Green function for the infinite square lattice is \cite{spitz}
\begin{equation}
G_{(p_1,q_1),(p_2,q_2)} \equiv G_{\vec{r}_1,\vec{r}_2} \equiv G(\vec{r}_2-\vec{r}_1) \equiv G_{0,0} + g_{p,q}, \quad \vec{r}_2-\vec{r}_1\equiv\vec{r}\equiv(p,q)
\end{equation}
with an irrelevant infinite constant $G_{0,0}$. The finite term $g_{p,q}$ is given explicitly by
\begin{equation}
g_{p,q}= \frac{1}{8\pi^2} \int\!\!\!\!\int_{-\pi}^{\pi}
\frac{e^{ \ci\, p \, \alpha + \ci \, q \, \beta}-1}{2-\cos\alpha-\cos\beta} \, \d \alpha \, \d \beta
\label{Green}
\end{equation}
and obeys the symmetry relations:
\begin{equation}
g_{p,q}=g_{q,p}=g_{-p,q}=g_{p,-q} \, .
\label{GreenSymmetry}
\end{equation}
After the integration over $\alpha$, it can be expressed in a more convenient form,
\begin{equation}
g_{p,q} = \frac{1}{4\pi} \int_{-\pi}^{\pi} \frac{t^p \, e^{ \ci\, q\, \beta} - 1 }{ \sqrt{y^2-1} } \,  \d \beta \, ,
\label{Green2}
\end{equation}
where $t = y - \sqrt{y^2-1}$, $y = 2 - \cos{\beta}$.

Below, we give $g_{p,q}$ for several values $p,q$ which are used in the text
\begin{equation}
\begin{array}{lllll}
  g_{0, 1} = - \frac{1}{4}    & \quad & g_{0, 2} = -1 + \frac{2}{\pi}            & \quad & g_{0, 3} = -\frac{17}{4} + \frac{12}{\pi} \\
  g_{1, 1} = - \frac{1}{\pi}  & \quad & g_{1, 2} = \frac{1}{4} - \frac{2}{\pi}   & \quad & g_{1, 3} = 2 - \frac{23}{3\pi} \\
  g_{2, 2} = - \frac{4}{3\pi} & \quad & g_{2, 3} = -\frac{1}{4} - \frac{2}{3\pi} & \quad & g_{3, 3} = - \frac{23}{15\pi} \, .
\end{array}
\label{someGF}
\end{equation}

\end{document}